\documentclass[12pt]{article}

\usepackage{amsfonts,amssymb,a4}

\newcommand{\qed}{\hfill\rule{3mm}{3mm}}
\newtheorem{theorem}{Theorem}[section]
\newtheorem{definition}{Definition}[section]
\newtheorem{proposition}{Proposition}[section]
\newtheorem{Example}{Example}[section]
\newenvironment{example}{\begin{Example}\rm }{\end{Example}}
\newtheorem{Remark}{Remark}[section]
\newenvironment{remark}{\begin{Remark}\rm }{\end{Remark}}

\makeatletter
\@addtoreset{equation}{section}
\makeatother

\newcommand{\SL}{{\rm sl}}
\newcommand{\real}{\mathbb{R}}
\newcommand{\thun}{\theta^1,\ldots,\theta^n}
\newcommand{\chds}{{\rm Ch}\bigl(D(\Sigma)\bigr)}
\newcommand{\om}[1]{\omega^{#1}}
\newcommand{\reff}[1]{(\ref{#1})}
\def\proof{\par\medskip\noindent{\sc Proof.\ }}
\newcommand{\isom}{\stackrel{\sim}{\to}{}}
\newcommand{\jj}{\mathcal J}
\newcommand{\mm}{\mathcal M}

\begin{document}

\begin{titlepage}
\title{Differential equations and moving frames}

\author{Odinete Ren\'ee Abib\\
{\small\it Laboratoire de Math\'ematiques Rapha\"el Salem}       \\[-1.7mm]
  {\small\it UMR 6085 CNRS}         \\[-1.7mm]
  {\small\it Universit\'e de Rouen}          \\[-1.7mm]
  {\small\it Avenue de l'Universit\'e, BP.12}          \\[-1.7mm]
  {\small\it 76801 Saint Etienne du Rouvray, FRANCE}      \\[-1.7mm]
  {\small\tt renee.abib@univ-rouen.fr}
}
\end{titlepage}

\maketitle

\begin{abstract}
We shall study the foundations of the differential geometric
consideration for differential equations.  We show a local structure theorem.
The main idea lies in the structure equations.  The Lie algebra
aspects of local differential equations is studied too.
\end{abstract}

\section{Introduction}

The purpose of the present paper is to study the relationship between
differential equations, Pfaffian systems and geometric structures, via
the method of moving frames of E.~Cartan~\cite{cartgroup,griffiths}

Following Cartan, we deal with every differential equation as a
Pfaffian system on a suitable manifold (Section~\ref{sec:3}).  This is
the fundamental idea of Cartan.  Further, we shall consider the
structure equations which are satisfied by Pfaffian systems determined
by differential equations.  The intergration of a given differential
equation is deeply related to the structure equation associated with
the differential equation.  We shall show it by means of some
examples.

In Section~\ref{sec:4}, we shall establish a local structure theorem
(Theorem~\ref{the:4.1}).  By virtue of this theorem, differential
equations can be regarded as a differential geometric structure on a
manifold.  In Section~\ref{sec:6}, we shall consider the Lie algebraic
aspect of local differential equations; each differential Lie algebra
(definition~\ref{def:6.1}) determines locally a local differential
equation (Theorem~\ref{the:6.1}); if $g$ is a semi-simple graded Lie
algebra, then $g$ has a structure of fundamental Lie algebra
(Theorem~\ref{the:6.2}).  Moreover $\SL(2,\real)$ has a structure of
differential Lie algebra which is not fundamental.  

In Section~\ref{sec:7} we study one system, which is one of the
typical examples in Cartan's paper~\cite{cartpfaff}, related with
$G$-structures and the local automorphism group of the given system.
The Section~\ref{sec:2} conscerns remarks on Pfaffian systems, Cauchy
characteristics and solvable systems.

In this paper, by the language differentiable we mean differentiable
of class $C^\infty$.  

I thank Marco Antonio Teixeira, Luiz San Martin, Paulo R\'egis Ruffino
for encouragement and IMECC-UNICAMP, BRASIL for their hospitality
during the preparation of this work.

\section{Cauchy characteristic system}\label{sec:2}

We begin with the preliminary remarks on Pfaffian systems.  Let $M$ be
a differentiable manifold.  $F(M)$ denotes the ring of real-valued
differentiable functions on $M$ and $\Lambda^1(M)$ the $F(M)$-module
of all 1-forms (Pfaffian forms) on $M$.  A $F(M)$-submodule $\Sigma$ of
$\Lambda^1(M)$ is said a \emph{Pfaffian system} of rank $n$ on $M$ if
$\Sigma$ is generated by $n$ linearly independent Pfaffian forms 
$\thun$.  A submanifold $N$ of $M$ is said an \emph{integral manifold}
of $\Sigma$ if $i^*\theta=0$ for all $\theta\in\Sigma$, where $i$
denotes the inclusion $N\hookrightarrow M$.  A differentiable function
$f$ on $M$ is said a \emph{first integral} of $\Sigma$ if the exterior
derivative $df$ belongs to $\Sigma$.   By the symbol
$\Sigma=\langle\thun\rangle$ we mean that the Pfaffian system $\Sigma$
is generated by the linearly independent Pfaffian forms $\thun$
defined on $M$. 

For each Pfaffian system $\Sigma$ on $M$, we can construct the dual
system, that is, the differentiable subbundle $D(\Sigma)$ of the
tangent bundle $T(M)$ on $M$ such that the fiber dimension of
$D(\Sigma)$ is equal to $\dim M-n$.  Let $\underline{D(\Sigma)}$ be
the sheaf of germs of local vector fields which belong to $D(\Sigma)$
and $\underline{D(\Sigma)}_{\,x}$ ($x\in M$) the stalk of 
$\underline{D(\Sigma)}$ at $x$.  We set
\[
\underline{\chds}_{\,x} \;=\; \Bigl\{
A\in \underline{D(\Sigma)}_{\,x} \,\Bigm| 
\Bigl[A\,,\,\underline{D(\Sigma)}_{\,x}
\Bigr]\subset \,\underline{D(\Sigma)}_{\,x}  \Bigr\}
\]
where $[\,,\,]$ denotes the natural bracket operation.  Further, for
each $x\in M$, we define the subspaces $\chds_x$ of $T_x(M)$ by
\[
\chds_x \;=\; \Bigl\{ X_x\in D(\Sigma)_x \Bigm| \, \underline X_{\,x}\in
\underline{\chds}_{\,x}\Bigr\}\;,
\]
where $X$ denotes a vector field and $\underline X_{\,x}$ the germ at
$x$ determined by $X$.  We suppose that $\dim \chds_x$ is constant on
$M$.  Thus, we obtain the subbundle $\chds$ of $T(M)$.  $\chds$ is
called the \emph{Cauchy charateristic} of $D(\Sigma)$.  The dual
system of $\chds$ is called the \emph{Cauchy characteristic system} of
$\Sigma$.  The following theorem is due to Cartan~\cite{cartinte,cartcomplet}.

\begin{theorem}\label{the:2.1}
Let $\Sigma=\langle\thun\rangle$ be a Pfaffian system.
\begin{enumerate}
\item If $\Sigma$ is completely integrable, i.e.\ $d\theta^i=0$ (mod.\
  $\thun$) $i=1,2,\ldots,n$, then ${\rm Ch}(\Sigma)=\Sigma$.
\item If $\Sigma$ is not completely integrable, then there exist
  linearly independent Pfaffian forms $\omega^1,\ldots,\omega^m$
  satisfying the following conditions:
\begin{itemize}
\item[(i)] $\thun, \omega^1,\ldots,\omega^m$ are also linearly
  independent; 
\item[(ii)] $(\thun, \omega^1,\ldots,\omega^m)$ forms a (local)
  generator of ${\rm Ch}(\Sigma)$;
\item[(iii)] $d\theta^i=\sum_{j,k=1}^m C_{jk}^i \omega^j\wedge\omega^k$
  (mod.\ $\thun$), where $C_{jk}^i$ denotes a differentiable function
  ($i=1,2,\ldots,n$; $j,k=1,2,\ldots,m$).
\end{itemize}

\item ${\rm Ch}(\Sigma)$ is completely integrable.

\item Let $x^1,\ldots, x^{n+m}$ be independent first integrals of ${\rm
    Ch}(\Sigma)$. Then there exist linearly independent Pfaffian forms
    $\overline\theta^i= \sum_{j=1}^{n+m} A_j^i(x^1,\ldots, x^{n+m})
    \,dx^j$, $i=1,2,\ldots,n$, such that $(\overline\theta^1,\ldots,
    \overline\theta^n)$ forms a (local) generator of $\Sigma$.

\end{enumerate}

\end{theorem}

By making use of property 2.(ii), we can construct the Cauchy
characteristic system ${\rm Ch}(\Sigma)$.

\begin{example} \label{exa:2.1} Consider the Pfaffian system
  $\Sigma=\langle\theta\rangle$, $\theta=dz+pdx+p^2dy$, on
  $\real^4=\bigl\{(x,y,z,p)\bigr\}$.  We have $dtheta=dp\wedge
  (dx+2p\,dy)$ and
\[ \omega^1=dp\quad,\quad \omega^2=dx+2p\,dy \quad,\quad \omega^3=p\quad;
\]
determine the Cauchy characteristic system of $\Sigma$.  We can find
by quadrature three independent first integrals as follows:
\[ u_1=z+xp+yp^2 \quad,\quad u_2=x+2yp \quad,\quad u_3=p\quad;
\]
and $\theta$ itself is expressed as $\theta=du_1-u_2\;du_3$.
\end{example}

\begin{definition}\label{def:2.1}
A system $(\omega^1,\ldots,\omega^m)$ of linearly independent Pfaffian
forms on $M$ will be said a \emph{solvable system} of
$\Sigma=\langle\thun\rangle$ if it satisfies the following conditions:
\begin{itemize}
\item[(i)] $(\omega^1,\ldots,\omega^m)$ forms a generator of ${\rm
    Ch}(\Sigma)$; 
\item[(ii)] $d\omega^1=0$ and $d\omega^p\equiv 0$ (mod.\
  $\omega^1,\ldots,\omega^{p-1}$) for all $p=2,3,\ldots,m$.
\end{itemize}
\end{definition}

If we can find a solvable system of $\Sigma$, then $m$ independent
first integrals of ${\rm Ch}(\Sigma)$ are given by quadrature.  In the
above example, the system $(\omega^1,\omega^2,\omega^3)$ is a solvable
system of $\Sigma=\langle\theta\rangle$.

\section{Differential equations and structure equations}\label{sec:3}

In this section we shall consider, by means of simple examples, the
relation between the differential equations and Pfaffian systems.

\paragraph{ a)} Take the first order equation on $\real^2=\{(x,y)\}$
\begin{equation}
\label{eq:3.1}
{\partial z\over \partial x} + \frac{1}{2} \Bigl(
{\partial z\over \partial y}\Bigr)^2 \;=\; 0\;.
\end{equation}
Setting on $\real^4=\{(x,y,z,q)\}$, $\om1=dx$, $\om2=dy-q\,dx$,
$\om3=dz+\frac12 q^2\,dy -q\,dy$, $\om4=dq$, we have
\begin{equation}
\label{eq:3.2}
\left\{\begin{array}{l}
d\om1 = 0\;,\\
d\om2 = \om1\wedge\om4\;,\\
d\om3 = \om2\wedge\om4\;,\\
d\om4 =0\;.
\end{array}\right.
\end{equation}
Each integral of \reff{eq:3.1} defines a 2-dimensional integral
manifold of $\langle\om3\rangle$ on which $\om1$ and $\om2$ are
linearly independent.  The equation \reff{eq:3.1} is left invariant by
the automorphism group of the absolute parallelism $\om1,\ldots,\om4$
on $\real^4$.  The structure of this group is determined by the
equation \reff{eq:3.2}.  The integration of the equation \reff{eq:3.1}
depends deeply on the structure equation \reff{eq:3.2} of this group.
In this case $(\om1,\om2,\om3)$ forms a solvable system of
$\langle\om3\rangle$.  Therefore three independent first integrals of
${\rm Ch}(\langle\om3\rangle)$ are given by quadrature as 
\[
u_1=q \quad,\quad u_2=z+\frac12 x q^2 - yq \quad,\quad
u_3=y-xq\quad;
\]
and we have $\om3=du_2 + u_3\,du_1$.  The formula
\[
\left\{\begin{array}{l}
z+\frac12xq^2 -yq = f(q)\;,\\
y-xq + f'(q) =0
\end{array}\right.
\]
gives an integral surface of the equation \reff{eq:3.1}, where $f$ is
a differentiable function and $f'$ denotes its derivative.

Conversely, we consider an absolute parallelism $\om1$, $\om2$,
$\om3$,  $\om4$ on $\real^4$ satisfying the equations
\begin{equation}
  \label{eq:3.3}
 \left\{\begin{array}{l}
d\om1 \equiv 0\;,\; d\om2\equiv 0\quad (\mbox{mod. } \om1,\om2)\;,\\
d\om3 \equiv \om2\wedge\om4 \quad (\mbox{mod. } \om3)\;.
 \end{array}\right.
\end{equation}
Let $x$ and $y$ be two independent first integals of the completely
integrable Pfaffian system $\om1=\om2=0$.  If we reduce $\om3$ to the
submanifold defined by the equations $x=$const., $y=$const., then from
the equation $d\om3\equiv\om2\wedge\om4$ (mod.\ $\om3$) we have
$d\om3\equiv 0$ (mod.\ $\om3$) on this submanifold.  Therefore $\om3$
must be of the form
\[
\om3\;=\; a(dz-p\,dx-q\,dy)\;,
\]
where $a$ is a non-zero function.  Since $\om1$, $\om2$, $\om3$ are
linearly independent, the functions $x$, $y$ and $z$ are also
independent. 

By this procedure we can determine the functions $p$ and $q$ of the
variables $x$, $y$, $z$ and another $t$:
\begin{equation}
  \label{eq:3.4}
  p=p(x,y,z,t)\quad,\quad q=q(x,y,z,t)\;;
\end{equation}
and the same equation $d\om3\equiv\om2\wedge\om4$ (mod.\ $\om3$)
implies 
\[
{\rm rank}\Bigl({\partial p\over \partial t}\,,\,
{\partial q\over \partial t}\Bigr)\;=\; 1\;.
\]
On a 2-dimensional integral manifold of $\om3=0$ on which $x$ and $y$
are still independent, $p$ and $q$ can be considered as the first
partial derivatives of $z=z(x,y)$.  Therefore the equation
\reff{eq:3.4} can be regarded as a first-order differential equation.
For example, the differential equation
\[
{\partial z\over \partial x} + \frac12\Bigl(
{\partial z\over \partial y} - f(x,y)\Bigr)^2
\;=\; g(x,y)
\]
belongs to the family determined by the structure equation
\reff{eq:3.2}, where $f(x,y)$ and $g(x,y)$ are differentiable
functions satisfying the equation
\[
{\partial f\over \partial x} \;=\;
{\partial g\over \partial y}\;.
\]

\paragraph{b)} Next, we consider an absolute parallelism $\om1$, 
$\om2$, $\om3$, $\om4$, $\om5$, $\om6$ on $\real^6$ satisfying
\begin{equation}
\label{eq:3.5}
\left\{\begin{array}{l}
d\om1 \equiv 0\;,\; d\om2\equiv 0\quad (\mbox{mod. } \om1,\om2)\;,\\
d\om3 \equiv \om1\wedge\om4+\om2\wedge\om5 \quad (\mbox{mod. } \om3)\;,\\
d\om4 =0\quad (\mbox{mod. } \om3,\om4,\om5)\;,\\
d\om5 = \om2\wedge\om6\quad (\mbox{mod. } \om3,\om4,\om5)\;.
\end{array}\right.
\end{equation}
Let $x$ and $y$ be two independent first integrals of the completely
Pfaffian integrable system $\om1=\om2=0$; $\om3$ is expressed as
\[
\om3\;=\; a(dz-p\,dx-qd\,y) \qquad (a\neq 0)\;.
\]
The functions $x$, $y$, $z$, $p$ and $q$ are independent first
integrals of the completely integrable Pfaffian system 
$\om1=\om2=\om3=\om4=\om5=0$.  Therefore $\om4$ and $\om5$ can be
written by means of the exterior derivatives $dx$, $dy$, $dz$, $dp$,
$dq$ and the formulas
\begin{eqnarray*}
  dp-r\,dx -s\,dy &=& a_1\om4 + a_2\om5 + a_3\om 3\\
dq -s'dx-t\, dy &=& a_4 \om4 + a_5\om5 + a_6\om3
\end{eqnarray*}
determine the functions $r$, $s$, $s'$, $t$ and $a_i$'s of the
variables $x$, $y$, $z$, $p$, $q$ and another $u$.  From the equation
$d\om3 \equiv \om1\wedge\om4+\om2\wedge\om5$ (mod.\ $\om3$), one can
verify that the function $s$ coincides with $s'$.  Moreover, the
equations $d\om4 =0\quad$, $d\om5 = \om2\wedge\om6$ (mod.\ 
$\om3,\om4,\om5$) imply
\[
{\rm rank}\Bigl({\partial r\over \partial u}\,,\,
{\partial s\over \partial u}\,,\,
{\partial t\over \partial u}\Bigr)\;=\; 1\;.
\]
Therefore the functions
\[
r=r(x,y,z,p,q,u) \;,\;
s=s(x,y,z,p,q,u) \;,\;
t=t(x,y,z,p,q,u) 
\] determine a system of second-order partial differential equations.
This family of systems of differential equations determined by an
absolute parallelism satisfying \reff{eq:3.5} is the main subject of
Cartn's researches in his paper \cite{cartpfaff}.  

For example, take the system of differential equations
(c.f.~\cite[\S\S\, 13, 14]{cartpfaff}) 
\begin{equation}
  \label{eq:3.6}
  {\partial^2 z\over\partial x^2}=0 \quad , \quad
{\partial^2 z\over \partial x\,\partial y} = z-x
{\partial z\over\partial x}\;.
\end{equation}
Putting on $\real^6=\{(x,y,z,p,q,t)\}$
$\om1=dx$, $\om2=dy$, $\om3=dz-p\,dx-q\,dy$,
$\om4=dp-(z-xp)\,dy$, $\om5=dq-(z-xp)\,dx-t\,dy$
and $\om6=dt-\bigl(q-x(z-xp)\bigr)\,dx$ we have the structure
equations 
\begin{eqnarray*}
\left\{\begin{array}{l}
d\om1 \equiv 0\;,\; d\om2\equiv 0\;,\\
d\om3 \equiv \om1\wedge\om4+\om2\wedge\om5 \;,\\
d\om4 \equiv \om2\wedge\om3-x\,\om2\wedge\om4 \;,\\
d\om5 = \om2\wedge\om6+ \om1\wedge\om3-x\,\om1\wedge\om4 \;,\\
d\om6 = \om1\wedge\om5-x\, \om1\wedge\om3-x^2\,\om1\wedge\om4 
+ K\,\om1\wedge\om2\;,
\end{array}\right.
\end{eqnarray*}
where $K=t-xq+x^2(z-xp)$.  The absolute parallelism satisfies the
equations \reff{eq:3.5}.  It is easy to see that the system
$(\om2,\om3,\om4,\om5,\om6)$ forms a solvable system of
$\Sigma=\langle\om3,\om4,\om5\rangle$.  Five independent first
integrals of the solvable system are given by quadrature as follows:
\[
u_1=y\quad,\quad u_2=z-xp\quad,\quad u_3=p
\quad,\quad u_4=q-x(z-xp)\quad,\quad u_5=K\quad,
\]
and we have (c.f.~\cite[\S 10, IV]{cartpfaff})
\begin{eqnarray*}
  \left\{\begin{array}{rcl}
\om3-x\,\om4 &=& du_2 -u_4\,du_1\,\\
\om4 &=& du_3 - u_2\,du_1\;,\\
\om5-x\,\om3 &=& du_4 -u_5\,du_1\;.
\end{array}\right.
\end{eqnarray*}
By this expression, the general integral surface of \reff{eq:3.6} is
given by the formulas:
\[
p=f(y) \;,\; z-xp=f'(y)\;,\; q-x(z-xp)=f''(y)\;,\;
t-x\bigl(q-x(z-xp)\bigr) = f'''(y)
\]
where $f$ is a differentiable function and $f'$, $f''$ and $f'''$
denote its derivatives.

\section{Differential geometric structures}\label{sec:4}

In the previous section we have seen that the integration of
differential equation is deeply related to the structure equations of
differential equations.  In this section we shall consider the
differential geometric structures for differential equations.

Let $V_{-1}$ and $V_0$ be finite dimensional real vector spaces.
We define by induction the real vector spaces $V_k$,
$k=1,2,\ldots$ as follows.  Let $V_1={\rm Hom}(V_{-1},V_0)$; 
$V_{k-1}$ ($k\ge 2$) being determined, we set
\[
V_k\;=\;\Bigl\{ X\in {\rm Hom}(V_{-1},V_{k-1})\, \Bigm|
X(u)(v)=X(v)(u)\,,\; u,v, \in V_{-1}\Bigr\}\;.
\]
We have $V_k\cong V_0\otimes S^k(V^*_{-1})$ as a vector space
($k=0,1,2,\ldots$), where $S^k(V^*_{-1})$ denotes the symmetric tensor
space of the dual space $V^*_{-1}$.  For an integer $k\ge 1$ we set
\[ W_k(V_{-1},V_0) \;=\; V_{-1}\oplus V_0\oplus \cdots\oplus V_k
\qquad \mbox{(direct sum),}
\]
and we define the bracket operation $[\,,\,]$ on $W_k(V_{-1},V_0)$ as
follows:
\begin{itemize}
\item[(i)] For all $X_{-1}\in V_{-1}$, $X_p\in V_p$ ($p\ge 1$),
\[ [X_p,X_{-1}]\;=\; -[X_{-1}, X_p] \;=\; X_p(X_{-1})\;; \]
\item[(ii)] $[X,Y]=0$ for any other combination
\end{itemize}
By this bracket operation, $ W_k(V_{-1},V_0)$ becomes a nilpotent Lie
algebra.  It is easy to prove the following.

\begin{proposition}\label{pro:4.1}\mbox{}
\begin{itemize}
\item[(i)] For a non-zero element $X_{-1}\in V_{-1}$,
  $[X_{-1},V_p]=V_{p-1}$ ($p\ge 1$).
\item[(ii)] If $[X_{-1},V_p]=(0)$ for $X_p\in V_p$ ($p\ge 1$), then
  $X_p=0$. 
\item[(iii)] For an arbitrary subspace $V_k^0$ of $V_k$,
\[ W_k^0(V_{-1},V_0) \;=\; V_{-1}\oplus V_0\oplus \cdots\oplus V_k^0
\qquad \mbox{(direct sum)}
\]
is a Lie subalgebra of $W_k(V_{-1},V_0)$.
\end{itemize}
\end{proposition}

\begin{example}\label{exa:4.1}\mbox{}

\begin{itemize}
\item[(1)] $\dim V_{-1}=1$, $\dim V_0=1$.  We have $\dim V_k=1$ for
any $k\ge 1$.  There exists a basis $X_{-1},X_0,X_1,\ldots,X_k$ of
$W_k(V_{-1},V_0)$ such that $X_p\in V_p$ ($-1\le p \le k$) and
$[X_{-1},X_p]=-X_{p-1}$ ($1\le p \le k$).

\item[(2)] $\dim V_{-1}=2$, $\dim V_0=1$.   We have $\dim V_k=k+1$ for
any $k\ge 1$.
\begin{itemize}
\item[(i)] $k=1$.  There exists a basis $X_1,X_2,X_3,X_4,X_5$ of 
$W_k(V_{-1},V_0)$ such that $X_1,X_2\in V_{-1}$; $X_3\in V_0$;
$X_4,X_5\in V_1$ and
\[ [X_1,X_4]=-X_3\quad,\quad [X_2,X_5]=-X_3 \]
and otherwise $[X,Y]=0$.
\item[(ii)] $k=2$.  There exists a basis
  $X_1,X_2,X_3,X_4,X_5,X_6,X_7,X_8$ of  
$W_2(V_{-1},V_0)$ such that $X_1,X_2\in V_{-1}$; $X_3\in V_0$;
$X_4,X_5\in V_1$ ; $X_6,X_7,X_8\in V_2$ and
\begin{eqnarray*}
&& [X_1,X_4]=-X_3\quad,\quad [X_2,X_5]=-X_3\quad,\quad 
 [X_1,X_6]=-X_4\\
&& [X_1,X_7]=-X_5 \quad,\quad
 [X_2,X_7]=-X_4\quad,\quad [X_2,X_8]=-X_5 
\end{eqnarray*}
and otherwise $[X,Y]=0$.
\end{itemize}
\end{itemize}
\end{example} 

Let $\pi:M\to N$ be a fibered manifold on a differentiable manifold
$N$ and $J^k(M,\pi)$ the space of $k$-jets of local sections of
$\pi$.  If $\dim N=\dim V_{-1}$ and $\dim M=\dim (V_{-1}\oplus
V_0)$, $W_k(V_{-1},V_0)$ is regarded as the local structure of
$J^k(M,\pi)$, i.e.\ $J^k(M,\pi)\cong W_k(V_{-1},V_0)$ (locally
diffeomorphic). 

Let $u^p:W_k(V_{-1},V_0)\to V_p$ ($-1\le p\le k$) be the natural
projection.  We regard $u^p$ as a vector-space valued function on
$W_k(V_{-1},V_0)$, so that the system $(u^{-1}, u^0,\ldots, u^k)$ can
be considered as a linear coordinate system on $W_k(V_{-1},V_0)$.  We
set $\theta^{-1}=du^{-1}$, $\theta^p=du^p-[u^{p-1},du^{-1}]$ ($0\le
p\le k-1$), $\theta^k=du^k$ and
$\theta=\theta^{-1}+\theta^0+\cdots +\theta^k$.  $\theta$ is a 
$W_k(V_{-1},V_0)$-valued 1-form on $W_k(V_{-1},V_0)$.  We have
\[
\left\{\begin{array}{l}
d\theta^{-1} = 0\\
d\theta^p + [\theta^{-1}\wedge \theta^{p+1}] = 0
\quad (0\le p \le k-1).
\end{array}\right.
\]
For example, making use of the notations in Example~\ref{exa:4.1},
2.(ii), we set $u^{-1}=x X_1+y X_2$, $u^0=zX_3$,
$u^1=pX_4+qX_5$, $u_2=r X_6+s X_7+tX_8$ and
$\theta=\sum_{i=1}^8 \om{i}X_i$.  Then we have 
$\om1=dx$, $\om2=dy$, $\om3=dz-p\,dx-q\,dy$,
$\om4=dp-r\,dx-s\,dy$, $\om5=dq-s\,dx-t\,dy$
$\om6=dr$, $\om7=ds$, $\om8=dt$ and
\[
\left\{\begin{array}{l}
d\om1 = 0\;,\; d\om2 = 0\;,\\
d\om3 = \om1\wedge\om4+\om2\wedge\om5 \;,\\
d\om4 = \om1\wedge\om6+\om2\wedge\om7 \;,\\
d\om5 = \om1\wedge\om7+ \om2\wedge\om8 \;.
\end{array}\right.
\]

Let $\rho^k:W_k(V_{-1},V_0)\to W_{k-1}(V_{-1},V_0)$
($k\ge 1$) be the natural projection, where we put 
$W_0(V_{-1},V_0)= V_{-1}\oplus V_0$; $W_k(V_{-1},V_0)$ can be
considered as a fibered manifold on $W_{k-1}(V_{-1},V_0)$ with the
fibering $\rho^k$. 

\begin{definition}\label{def:4.1}
We shall say that a submanifold $R_k$ of $W_k(V_{-1},V_0)$ is a
\emph{local differential equation} of order $k$ if $R_k$ admits an
absolute parallelism and if there exist an open submanifold $U$ of 
$W_{k-1}(V_{-1},V_0)$ such that $\rho^k(R_k)=U$ and
$\rho^k\bigr|_{R_k} : R_k\to U$ is a fibered submanifold of
$\rho^k:(\rho^k)^{-1}(U)\to U$.
\end{definition}

Let $i: R_k\hookrightarrow W_k(V_{-1},V_0)$ be the inclusion and 
$\omega=i^*\theta$ the induced $W_k(V_{-1},V_0)$-valued 1-form on
$R_k$.  According to the direct sum decomposition of
$W_k(V_{-1},V_0)$, we decompose $\omega$ as
$\omega=\om{-1}+\om0+\cdots
\widetilde\omega^k$, where $\om{p}$ (resp.\ $\widetilde\omega^k$) is
a $V_p$-valued 1-form (resp.\ $V_k$-valued 1-form) on $R_k$ ($-1\le
p\le k-1$).  Let $n$ be the fiber dimension of $\rho:R_k\to U$.  Then
there exist $n$ linearly independent Pfaffian forms
$\omega_1^k,\ldots,\omega_n^k$ which are also linearly independent of
the Pfaffian forms obtained from $\om{-1},\om0,\ldots,\om{k-1}$.  We
fix a $n$-dimensional vector subspace $V_k^0$ of $V_k$ and its basis
$X_1,\ldots,X_n$ and we set $\om{k}=\sum_{j=1}^n \omega_j^k\,X_j$.
$\om{k}$ is a $V_k^0$-valued 1-form on $R_k$.  We define the
differentiable mappings
\begin{eqnarray*}
  F_p^k &:& R_k \longrightarrow {\rm Hom}(V_p,V_k) \qquad (-1\le p\le
  k-1)\;,\\ 
F_k^k &:& R_k \longrightarrow {\rm Hom}(V_k^0,V_k)
\end{eqnarray*}
by the formula
\[ d(u^k\circ i) \;=\; F_{-1}^k(\om{-1}) +\cdots+
F_k^k(\om{k})\;;
\]
and we define the differentiable mapping
\[ T: R_k \longrightarrow {\rm Hom}(V_{-1}\times V_k^0,V_{k-1}) 
\]
by the formula 
\[ T_x(X_{-1},X_k) \;=\; [X_{-1}, F_k^k(x)(X_k)]
\qquad (x\in R_k\,,\, X_{-1}\in V_{-1}\,,\, X_k\in V_k^0)\;.
\]
Since the rank of the inclusion $i$ is maximal on $R_k$, the linear
mapping $F_k^k(x):V_k^0\to V_k$ ($x\in R_k$) is injective.  Therefore
$T$ has the following property:
\begin{itemize}
\item[(C$_1$)] For each $x\in R_k$, $T_x(V_{-1},X_k) =0$ ($X_k\in
  V_k^0$) implies $X_k=0$.
\end{itemize}
It is also easy to prove the following properties:
\begin{itemize}
\item[(C$_2$)]
\begin{itemize}
\item[(i)] $d\om{-1}\equiv 0$ (mod.\ $\om{-1}$);
\item[(ii)] ($k\ge 2$) For $p=0,1,\ldots,k-2$,
\[ d\om{p} + [\om{-1}\wedge \om{p-1}] \;\equiv\; 0 \qquad 
(\mbox{mod. } \om0,\ldots,\om{p})\;;
\]
\item[(iii)] $d\om{k-1} + T(\om{-1}\wedge \om{k}) \equiv 0$ (mod.\ 
$\om{-1}\wedge \om{-1},\om0,\ldots,\om{k-1}$).
\end{itemize}
\end{itemize}

We have thus proved that for each local differential equation of order
$k$ there exist a differentiable mapping 
$T: R_k \to {\rm Hom}(V_{-1}\times V_k^0:V_{k})$ and an absolute
parallelism $\omega=\om{-1}+\om0+\cdots\om{k}$ satisfying the above
conditions (C$_1$), (C$_2$).

\begin{theorem}\label{the:4.1}
Let $V_k^0$ be a subspace of $V_k$ and put 
$W_k^0=V_{-1}\oplus V_0\oplus\cdots\oplus V_k^0$ (direct sum).  Let $R_k$
be a differentiable manifold with $\dim R_k = \dim W_k^0$.  If there
exists a $W_k^0$-valued absolute parallelism
$\omega=\om{-1}+\om0+\cdots\om{k}$ and a differentiable mapping 
$T: R_k \to {\rm Hom}(V_{-1}\times V_k^0:V_{k})$ satisfying the
conditions (C$_1$) and (C$_2$), then $R_k$ can be locally embedded
into $W_k(V_{-1},V_0)$ as a local differential equation or order $k$.  
\end{theorem}

\proof
Since $\omega$ gives rise to an isomorphism $\omega_x:T_x(R_k)\isom W_k^0$ 
($x\in R_k$) one can consider the inverse mapping of $\omega_x$, say
$\tau_x:W_k^0\isom T_x(R_k)$.  $\tau$ has the property:
$\om{p}\bigl((\tau(X_q)\bigr) = \delta_q^p X_q$, $X_q\in V_q$
($-1\le p,q\le k$).  By the condition (i) of C$_2$, we can find
differentiable mappings $v^{-1}: R_k\to V_{-1}$ and
$A_{-1}: R_k\to {\rm GL}(V_{-1})$  such that $\om{-1}=A_{-1}(dv^{-1})$.
Since the system $\om{-1}=\om0=0$ is completely integrable, there exists a
differentiable mapping $v^0:R_k\to V_0$ such that the system
$dv^{-1}=dv^0=0$ is equivalent to the system $\om{-1}=\om0=0$.
Therefore, $\om0$ can be written as
\[ \om0\;=\; A_0(dv^0-v^1\,dv^{-1})\;, \]
where $A_0$ denotes a differentiable mapping $A_0:R_k\to {\rm GL}(V_0)$
and $v^1$ denotes a differentiable mapping $v^1: R_k\to {\rm Hom}(V_{-1},V_0)=V_1$.
If $k=1$ the argument comes to an end.  Let $k\ge 2$.  Consider the following
proposition (P$_j$) for $1\le j\le k$:
\begin{itemize}
\item[(P$_j$)] There exist differentiable mappings $v^{-p}:R_k\to V_p$ and
$A_p:R_k\to {\rm GL}(V_p)$, $p=-1,0,1,\ldots j$, such that
\begin{itemize}
\item[(1)] $dv^{-1}, dv^0,\ldots,dv^{j-1}$ are linearly independent;
\item[(2)] $\om{-1}=A_{-1}(dv^{-1})$, $\om0\;=\; A_0(dv^0-v^1\,dv^{-1})$ and
for $p=1,2,\ldots,j-1$, $\om{p}\;=\; A_p(dv^p-v^{p-1}\,dv^{-1})$ 
(mod.\ $\om0,\ldots,\om{p-1}$).
\end{itemize}
\end{itemize}
We have proved (P$_1$).  For an integer $1\le j\le k-1$, assume that
(P$_j$) is  established.  From the inequality $0\le j-1\le k-2$ and
condition (ii) of (C$_2$) we have
\[ d\om{j-1}\;\equiv\; -[\om{-1}\wedge\om{j}] \qquad (\mbox{mod.\ }
\om0,\ldots,\om{j-1})\;.
\]
From (2) of (P$_j$) we have
\[ d\om{j-1}\;\equiv\; -A_{j-1}(dv{j}\wedge dv^{-1}) \qquad (\mbox{mod.\ }
\om0,\ldots,\om{j-1})\;.
\]
These two equations yield
\begin{equation}
\label{eq:4.1}
-A_{j-1}(dv{j}\wedge dv^{-1})\;\equiv\; -[\om{-1}\wedge\om{j}] \qquad (\mbox{mod.\ }
\om0,\ldots,\om{j-1})\;.
\end{equation}
Substituting $\tau(X_p)\wedge\tau(X_{-1})$, $X_{-1}\in V_{-1}$, $X_p\in V_p$
($j+1\le p\le k$) to this equation, we obtain
\[
\left(dv^j\bigl(\tau(X_p)\bigr)\right) \left(A_{-1}^{-1}(X_{-1})\right) \;=\; 0
\]
and hence
\begin{equation}
\label{eq:4.2}
dv^j\bigl(\tau(X_p)\bigr) \;=\; 0 \qquad (X_p\in V_p : p=j+1,\ldots,k)\;.
\end{equation}
Substituting $\tau(X_j)\wedge\tau(X_{-1})$, $X_{-1}\in V_{-1}$, $X_j\in V_j$
to equation \reff{eq:4.1}, we have 
\[
A_{j-1}\left(dv^j\bigl(\tau(X_j)\bigr)\right) \left(A_{-1}^{-1}(X_{-1})\right) \;=\; 
-[X_{-1},X_j]
\]
and hence
\begin{equation}
\label{eq:4.3}
dv^j\bigl(\tau(X_j)\bigr) \;=\; A_{j-1}^{-1}\circ X_j\circ A_{-1}\;.
\end{equation}
The equation \reff{eq:4.2} implies that $dv^j$ is expressed as
\begin{equation}
\label{eq:4.4}
dv^j\ \;=\; \; B_j(\om{j}) \qquad (\mbox{mod.\ } \om{-1},\om0,\ldots,\om{j-1})
\end{equation}
and the equation \reff{eq:4.3} implies that the differentiable mapping
$B_j:R_k\to {\rm Hom}(V_j,V_j)$ is given by the formula
\[
B_j(X_j) \;=\; A_{j-1}^{-1}\circ X_j\circ A_{-1} \qquad (X_j\in V_j)\;,
\]
so that $B_j(x)$ is non-singular for any $x\in R_k$.  Therefore one can see
that $dv^{-1},dv^0,\ldots,dv^j$ are linearly independent.  By equation \reff{eq:4.4},
$\om{j}$ can be written as
\[
\om{j}\;=\; A_j(dv^j-v^{j+1}\,dv^{-1}) \qquad \mbox{mod.\ } \om0,\ldots,\om{j-1})\;.
\]
where $A_j=B_j^{-1}$ and $v^{j+1}$ denotes a differentiable  mapping
$v^{j+1}:R_k\to {\rm Hom}(V_{-1},V_j)$.
Substituting $\tau(X_{-1})$, $X_{-1}\in V_{-1}$, to this equation, we have
\[
dv^j\bigl(\tau(X_{-1})\bigr) \;=\; v^{j+1}\bigl(A_{-1}^{-1}(X_{-1})\bigr) \;.
\]
Substituting $\tau(X_{-1})\wedge\tau(Y_{-1})$, $X_{-1},Y_{-1}\in V_{-1}$, 
to equation \reff{eq:4.1}, we obtain 
\[
dv^j\bigl(\tau(X_{-1})\bigr) \bigl(A_{-1}^{-1}(Y_{-1})\bigr) \;=\;
dv^j\bigl(\tau(Y_{-1})\bigr) \bigl(A_{-1}^{-1}(X_{-1})\bigr)\;.
\]
These two equations imply
\[
v^{j+1}(X_{-1})(Y_{-1}) \;=\; v^{j+1}(Y_{-1})(X_{-1})
\]
for any $X_{-1},Y_{-1}\in V_{-1}$, so that $v^{j+1}(x)$ lies in $V_{j+1}$ for any
$x\in R_k$. Thus we can establish by induction the proposition (P$_{j+1}$)
and hence (P$_k$).  Define the differentiable mapping $F:R_k\to W_k(V_{-1},V_0)$
by the formula 
\[
u^p\circ F \;=\; v^p \qquad (-1\le p \le k)
\]
and put
\[
dv^k\;\equiv\; v^k_k(\om{k}) \qquad (\mbox{mod.\ }\om{-1},\om0,\ldots,\om{k-1})\;.
\]
If $v_k^k(x)\in {\rm Hom}(V_k^0,V_k)$ is injective for any $x\in R_k$, $F$ is an immersion
and determines locally an embedding.  From the proposition (P$_k$) we have
\begin{eqnarray*}
d\om{k-1} &\equiv& -A_{k-1}\bigl(dv^k\wedge A_{-1}^{-1}(\om{-1})\bigr)
\qquad (\mbox{mod.\ }\om0,\ldots,\om{k-1})\\
&\equiv& -A_{k-1}\bigl(v_k^k(\om{k})\wedge A_{-1}^{-1}(\om{-1})\bigr)
\qquad (\mbox{mod.\ }\om{-1}\wedge\om{-1},\om0,\ldots,\om{k-1})\;.
\end{eqnarray*}
By the condition (iii) of (C$_2$) we obtain
\[
A_{k-1}\bigl(v_k^k(\om{k})\wedge A_{-1}^{-1}(\om{-1})\bigr)\;\equiv\;
T(\om{-1}\wedge\om{k})
\qquad (\mbox{mod.\ }\om{-1}\wedge\om{-1},\om0,\ldots,\om{k-1})\;.
\]
Substituting $\tau(X_{-1})\wedge\tau(X_k)$, $X_{-1}\in V_{-1}$, $X_k\in V_k^0$
to this equation, we have
\[
A_{k-1}\bigl(v_k^k(X_k)\wedge  A_{-1}^{-1}(X_{-1})\bigr)\;=\; T(X_{-1},X_k)\;.
\]
If $v_k^k(X_k)=0$, then $T(X_{-1},X_k)=0$ for any $X_{-1}\in V_{-1}$.
From condition (C$_1$) we obtain $X_k=0$.  Hence $v_k^k(x)\in{\rm Hom}(V_k^0,V_k)$
is injective for any $x\in R_k$.  Set 
$\widetilde\rho=\rho^k\circ F$.  By the definition of $F$
we have
\[
u^p\circ\widetilde\rho \;=\; u^p\circ\rho^k\circ F\;=\; u^p\circ F \;=\; v^p
\qquad (-1\le p\le k-1)\;.
\]
This relation and (1) of the proposition (P$_k$) imply that 
$\widetilde\rho$ is a submersion.  Therefore $F$ determines locally a local
differential equation of order $k$. \qed

\begin{remark}
By virtue of this theorem, a system $(R_k,W_k^0,T,\omega)$
satisfying the conditions stated in the theorem may be also called
a local differential equation of order $k$.
\end{remark}

\section{Equivalence}\label{sec:5}
Let $V_k^0$ be a subspace of $V_k$ ($k\ge 1$).  We set
$W_k^0=V_{-1}\oplus V_0\oplus\cdots\oplus V_k^0$
and $D^p=V_p\oplus V_{p-1}\oplus\cdots\oplus V_k^0$,
$p=0,1,\ldots, k$.  We define the Lie subgroup ${\rm G}(W_k^0)$
of ${\rm GL}\bigl(W_k(V_{-1},V_0)\bigr)$ as follows:
\[
{\rm G}(W_k^0) \;=\; \Bigl\{ g\in {\rm GL}\bigl(W_k(V_{-1},V_0)\bigr)
\Bigm| g(V_{-1}\oplus V_k^0)=V_{-1}\oplus V_k^0\;,\;
g(D^p)=D^p\; (0\le p\le k)\Bigr\}\;.
\]

\begin{definition}\label{def:5.1}
We shall say that two local differential equations $(R_k,W_k^0,T,\omega)$
and $(R_k,W_k^0,T',\omega')$ are \emph{structurally equivalent} if
there exists a differentiable mapping $A:R_k\to {\rm G}(W_k^0)$
such that $\omega'=A(\omega)$.
\end{definition}

\begin{definition}\label{def:5.2}
A local differential equations $(R_k,W_k^0,T,\omega)$ will be said
\emph{of type} $W_k^0$ if $T_x(X_{-1},X_k)=[X_{-1},X_k]_0$
for all $x\in R_k$, $X_{-1}\in V_{-1}$, $X_k\in V_k^0$ and if 
$\omega$ satisfies the condition:
\begin{itemize}
\item[(C$_2^{\,\prime}$)]
\begin{itemize}
\item[(i)] $d\om{-1}\equiv 0$ (mod.\ $\om{-1}$);
\item[(ii)] For $p=0,1,\ldots,k-1$,
\[ d\om{p} + [\om{-1}\wedge \om{p-1}] \;\equiv\; 0 \qquad 
(\mbox{mod. } \om0,\ldots,\om{p})\;,
\]
\end{itemize}
\end{itemize}
where $[\cdot,\cdot]_0$ denotes the natural bracket operation $W_k^0$
(cf.\ Proposition \ref{pro:4.1}).
\end{definition}

\begin{remark}\label{rem:5.1}
Almost all local differential equations which admit a lot of solutions
turn out to be structurally equivalent to a local differential equation
of type $W_k^0$ for some $V_k^0$.
\end{remark}

\begin{example}\label{exa:5.1}
Take the system of second order differential equations
\begin{equation}
  \label{eq:5.1}
  {\partial^2 z\over\partial x^2}=0 \quad , \quad
{\partial^2 z\over \partial x\,\partial y} = z\;.
\end{equation}
Putting on $\real^6=\{(x,y,z,p,q,t)\}$
$\om1=dx$, $\om2=dy$, $\om3=dz-p\,dx-q\,dy$,
$\om4=dp-z\,dy$, $\om5=dq-z\,dx-t\,dy$
and $\om6=dt$, we have 
\[
\left\{\begin{array}{l}
d\om1 = 0\;,\; d\om2\equiv 0\;,\\
d\om3 = \om1\wedge\om4+\om2\wedge\om5 \;,\\
d\om4 = \om2\wedge\om3+p\,\om1\wedge\om2 \;,\\
d\om5 = \om2\wedge\om6- \om1\wedge\om3-q\,\om1\wedge\om2 \;,\\
d\om6 = 0\;.
\end{array}\right.
\]
Since one can not remove the terms $p\om1\wedge\om2$ and $q\om1\wedge\om2$, 
$\omega$ is not of type $W_2^0$ for any $V_2^0\subset V_2$.  On the other hand,
the given system \reff{eq:5.1} has no solutions except $z=0$.
\end{example}

\begin{proposition}\label{pro:5.1}
Let $P$ be a differentiable manifold with $\dim P\ge \dim W_k^0$.
Suppose there exists a $W_k^0$-valued 1-form 
$\omega=\om{-1}+\om0+\cdots+\om{k}$ such that 
$\omega_p:T_p(P)\to W_k^0$ is surjective for any $p\in P$ and
$d\omega\equiv 0$ (mod.\ $\omega$).  If $\omega$
satisfies the condition (C$_2^{\,\prime}$) in Definition~\ref{def:5.2},
then $(P,W_k^0,\omega)$ determines locally a local differential
equation of type $W_k^0$.
\end{proposition}

\proof
Since $\omega=0$ is completely integrable, there exists, for each $p\in P$,
an open neighborhood $U$ of $p$, a differentiable manifold $R_k$ with
$\dim R_k=\dim W_k^0$ and a fibering $\pi: U\to R_k$ such that each
fiber is a maximal integral manifold of $\omega\bigr|_{\textstyle U}=0$.
Let $\sigma:R_k\to U$ be a differentiable cross section of $\pi$ 
and put $\overline\omega=\sigma^*\omega$.  Then it is clear that
$(R_k,W_k^0,\overline\omega)$ is a local differential equation of
type $W_k^0$.  In general, the obtained system depends on the
choice of cross sections. \qed

For a subspace $V_k^0$ of $V_k$ we set
\[
(V_k^0)^{(1)} \;=\; \Bigl\{ X\in {\rm Hom}(V_{-1},V_k^0) \Bigm|
X(u)(v)=X(v)(u)\;,\; u,v\in V_{-1}\Bigr\}\;.
\]
For a subspace $U$ of $V_{-1}$, we set
\[
V_k^0(U)\;=\; \Bigr\{ X_k\in V_k^0 \Bigm| X_k(u)=0\;,\; u\in U \Bigr\}\;.
\]

\begin{definition}\label{def:5.3}
A subspace $V_k^0$ of $V_k$ is said \emph{involutive} if there exists
a series of subspaces $(0)\subset U_0\subset U_1\subset\cdots
\subset U_{n-1}\subset U_n=V_{-1}$ with $\dim U_i=i$ such that
\[
\dim (V_k^0)^{(1)} \;=\; \sum_{i=0}^n \dim V_k^0(U_i)\;.
\]
A local differential equation $(R_k,W_k^0,\omega)$ of type $W_k^0$
is said \emph{involutive} if $V_k^0$ is involutive.
\end{definition}

\begin{example}\label{exa:5.2}
We use the notation of Example \ref{exa:4.1}.  By the symbol
$U=(X_1,X_2,\ldots,X_n)$  we mean that the vector space $U$
is spanned by the basis $X_1,X_2,\ldots,X_n$.
\begin{itemize}
\item[(1)] $\dim V_{-1}=1$, $\dim V_0=1$. 
\begin{itemize}
\item[(i)] $k=1$. Then $W_1^0=V_{-1}\oplus V_0=(X_1,X_2)$ with
$[X_1,X_2]=0$.  A $W_1^0$-valued 1-form $\omega=\om1X_1+\om2X_2$
is of type $W_1^0$ if it satisfies 
\[
\left\{\begin{array}{l}
d\om1 \equiv 0 \quad (\mbox{mod.\ } \om1)\;,\\
d\om2 \equiv 0 \quad (\mbox{mod.\ } \om2)\;.
\end{array}\right.
\]
\item[(ii)] $k=2$. Then $W_2^0=V_{-1}\oplus V_0\oplus V_1 =(X_1,X_2,X_3)$ with
$[X_1,X_3]=-X_2$ and otherwise $[X_i,X_j]=0$.  $\omega=\om1X_1+\om2X_2+\om3X_3$
is of type $W_2^0$ if it satisfies 
\[
\left\{\begin{array}{l}
d\om1 \equiv 0 \quad (\mbox{mod.\ } \om1)\;,\\
d\om2 = \om1\wedge\om3 \quad (\mbox{mod.\ } \om2)\;\\
d\om3 = 0 \quad (\mbox{mod.\ } \om2,\om3)\\\;.
\end{array}\right.
\]
The Pfaffian system $\Sigma=\langle\om2,\om3\rangle$ determines a family
of second-ordeer ordinary differential equations.
\end{itemize}
\item[(2)] $\dim V_{-1}=2$, $\dim V_0=1$. 
\begin{itemize}
\item[(i)] $k=1$. Let $V_1^0$ be a 1-dimensional subspace of $V_1$.

The we can choose a basis $X_1,X_2,X_3,X_4$ of  
$W_1^0$ such that $X_1,X_2\in V_{-1}$ and $[X_2,X_4]=-X_3$
and otherwise $[X_i,X_j]=0$.  $\omega=\sum_{i=1}^4 \om{i}X_i$
is of type $W_1^0$ if it satisfies
\[
\left\{\begin{array}{l}
d\om1 \equiv 0 \;,\; d\om2 \equiv 0\qquad (\mbox{mod.\ } \om1,\om2)\;,\\
d\om3 \equiv \om2\wedge\om4 \quad (\mbox{mod.\ } \om3)\;.
\end{array}\right.
\]
We have already seen this structure equation in Section~\ref{sec:3}.

\item[(ii)] $k=2$. Let $V_2^0$ be a 2-dimensional subspace of $V_2$.
Then the Lie algebra $W_2^0$ is isomorphic to the following three
Lie algebras.
\begin{itemize}
\item[(a)] $W_2^0=(X_1,X_2,X_3,X_4,X_5,X_6,X_7)$ with
\begin{eqnarray*}
&& [X_1,X_4]=-X_3\quad,\quad [X_2,X_5]=-X_3\quad,\quad 
 [X_1,X_6]=-X_4\\
&& [X_1,X_7]=-X_5 \quad,\quad
 [X_2,X_7]=-X_4\quad,\quad [X_2,X_6]=X_5 
\end{eqnarray*}
and otherwise $[X_i,X_j]=0$; $\omega=\sum_{i=1}^7 \om{i}X_i$
is of type $W_2^0$ if it satisfies
\[
\left\{\begin{array}{l}
d\om1 \equiv 0 \;,\; d\om2 \equiv 0\qquad (\mbox{mod.\ } \om1,\om2)\;,\\
d\om3 \equiv \om1\wedge\om4 +\om2\wedge\om5\quad (\mbox{mod.\ } \om3)\;,\\
d\om4 \equiv \om1\wedge\om6 +\om2\wedge\om7\quad (\mbox{mod.\ } \om3,\om4,\om5)\;,\\
d\om5 \equiv \om1\wedge\om7 -\om2\wedge\om6\quad (\mbox{mod.\ } \om3,\om4,\om5)\;.
\end{array}\right.
\]

\item[(b)] $W_2^0=(X_1,X_2,X_3,X_4,X_5,X_6,X_7)$ with
\begin{eqnarray*}
&& [X_1,X_4]=-X_3\quad,\quad [X_2,X_5]=-X_3\\
 &&[X_1,X_6]=-X_4\quad,\quad 
 [X_2,X_7]=-X_5 
\end{eqnarray*}
and otherwise $[X_i,X_j]=0$; $\omega=\sum_{i=1}^7 \om{i}X_i$
is of type $W_2^0$ if it satisfies
\[
\left\{\begin{array}{l}
d\om1 \equiv 0 \;,\; d\om2 \equiv 0\qquad (\mbox{mod.\ } \om1,\om2)\;,\\
d\om3 \equiv \om1\wedge\om4 +\om2\wedge\om5\quad (\mbox{mod.\ } \om3)\;,\\
d\om4 \equiv \om1\wedge\om6 \quad (\mbox{mod.\ } \om3,\om4,\om5)\;,\\
d\om5 \equiv \om2\wedge\om7\quad (\mbox{mod.\ } \om3,\om4,\om5)\;.
\end{array}\right.
\]
Let $V_2^0$ be a 1-dimensional involutive subspace of $V_2$.  Then there is only
one case up to isomorphic algebra.

\item[(c)] $W_2^0=(X_1,X_2,X_3,X_4,X_5,X_6)$ with
\[
 [X_1,X_4]=-X_3\quad,\quad [X_2,X_5]=-X_3
 \quad,\quad  [X_2,X_6]=-X_5 
\]
and otherwise $[X_i,X_j]=0$; $\omega=\sum_{i=1}^6 \om{i}X_i$
is of type $W_2^0$ if it satisfies
\[
\left\{\begin{array}{l}
d\om1 \equiv 0 \;,\; d\om2 \equiv 0\qquad (\mbox{mod.\ } \om1,\om2)\;,\\
d\om3 \equiv \om1\wedge\om4 +\om2\wedge\om5\quad (\mbox{mod.\ } \om3)\;,\\
d\om4 \equiv 0 \quad (\mbox{mod.\ } \om3,\om4,\om5)\;,\\
d\om5 \equiv \om2\wedge\om6\quad (\mbox{mod.\ } \om3,\om4,\om5)\;.
\end{array}\right.
\]
We have already discussed this case in Section~\ref{sec:3}.
\end{itemize}
\end{itemize}
\end{itemize}

\end{example}

\section{Lie algebraic aspects of differential equations}\label{sec:6}

In this section we shall consider the Lie-algebraic aspect of local
differential equations.  Let $V_k^0$ be a subspace of $V_k$.  We set
\[
W_k^0(V_k^0) = V_{-1}\oplus V_0\oplus \cdots\oplus V_k^0
\quad,\quad 
D^p = V_p\oplus V_{p-1}\oplus \cdots\oplus V_k^0
\qquad (0\le p\le k)\;.
\]
We define the Lie algebra $\jj\bigl(W_k^0(V_k^0)\bigr)$ as follows:
\[
\jj_k\bigl(W_k^0(V_k^0)\bigr) \;=\; \Bigl\{ X\in {\rm gl}\bigl(W_k^0(V_k^0)\bigr)
\Bigm| X(V_{-1}\oplus V_k^0) \subset V_{-1}\oplus V_k^0\;,\;
X(D^p)\subset D^p\; (0\le p\le k) \Bigr\}\;.
\]
Every element of $\jj_k\bigl(W_k^0(V_k^0)\bigr)$ is of the form:
\begin{eqnarray*}
&&\quad
V_{-1} \, V_0 \,V_1 \,\cdots \,V_{k-1} \,V_k^0
\\
&&\left(\begin{array}{cccccc}
*&0&0&\cdots &0&0\\
0&*&0&\cdots &0&0\\
0&*&*&\cdots&0&0\\
\vdots &\vdots&\vdots&\ddots&&\vdots\\
0&*&*&\cdots &*&0\\
*&*&*&\cdots &*&*
\end{array}
\right)\hspace{-.1cm}
\begin{array}{c}
V_{-1} \\ V_0 \\V_1 \\ \vdots \\V_{k-1} \\V_k^0
\end{array}
\end{eqnarray*}
where $*$ denotes a certain matrix. We define a mapping
$\partial: {\rm Hom}\bigl(W_k^0, \jj_k(W_k^0)\bigr) \to
{\rm Hom}(W_k^0,\wedge W_k^0:W_k^0)$ by the formula
\[(\partial S)(X\wedge Y) \;=\; S(X)(Y) - S(Y)(X)\;,
\]
for all $S\in {\rm Hom}\bigl(W_k^0, \jj_k(W_k^0)\bigr)$ and all $X,Y\in W_k^0$,
where we put $W_k^0=W_k^0(V_k^0)$.
\begin{definition}\label{def:6.1}
Let $\jj$ be a Lie algebra.  A system $(\jj,\mm,\jj_0)$ will be called a 
\emph{differential Lie algebra} if it satisfies the following conditions:
\begin{itemize}
\item[(1)] $\jj_0$ is a Lie subalgebra of $\jj$ and $\mm$ is a vector 
subspace of $\jj$ such that $\jj=\mm\oplus \jj_0$ (direct sum);
\item[(2)] For some subspaces $V_{-1}$, $V_0$ and $V_k^0$ of $\mm$, we have
$\mm=W_k^0(V_k^0)$;
\item[(3)] The linear isotropy representation $\rho:\jj_0\to {\rm gl}(\mm)$
($\rho(X_0)(X)=$ the $m$-component of $[X_0,X]$, $X_0\in\jj_0$, $X\in\mm$)
has its values in $\jj_k\bigl(W_k^0(V_k^0)\bigr)$;
\item[(4)] Let $\alpha:\mm\wedge\mm \to \mm$ be the linear mapping defined by
the formula 
\[
\alpha(X\wedge Y) \;=\;  \mbox{the $m$-component of } [X,Y]
\quad,\quad X,Y\in\mm\;.
\]
There exists an element $S\in {\rm Hom}\bigl(\mm,\jj_k(W_k^0)\bigr)$ such that
\[
\alpha(X\wedge Y) \;=\; [X,Y]_0+(\partial S)(X\wedge Y)\;,
\]
where $[\,,\,]_0$ denotes the natural bracket operation on $W_k^0(V_k^0)$.
If $J_0=(0)$, $\jj$ ($=\mm$) will be called fundamental.
\end{itemize}
\end{definition}

\begin{theorem}\label{the:6.1}
Let $(\jj,\mm,\jj_0)$ be a differential Lie algebra.  Let $G$ be a Lie group
corresponding to the Lie algebra $\jj$ and let $\theta$ be the
Maurer-Cartan form of $G$.  Then the $\mm$ component $\omega$ of 
$\theta$ with respect to the decomposition $\jj=\mm\oplus\jj_0$ determines
locally a local differential equation.
\end{theorem}

\proof
Let $\mm=W_k^0(V_k^0)$ for the subspaces $V_{-1}$, $V_0$ and $V_k^0$ of $\mm$.
We set $\theta=\omega+\omega_0$, $\omega_0$ being the $\jj_0$-component of $\theta$.
Then from the structure equation of Maurer-Cartan we obtain
\[
d\omega +\frac12\,\alpha(\omega\wedge\omega) - \rho(\omega_0)\wedge\omega
\;=\;0
\]
and hence
\[
d\omega+\frac12[\omega\wedge\omega]_0 +
\frac12(\partial S)(\omega\wedge\omega) - 
\rho(\omega_0)\wedge\omega
\;=\;0\;.
\]
Since $S(X)$ and $\rho(X_0)$ ($X\in\mm$, $X_0\in\jj_0$) lie in 
$\jj_k(W_k^0)$, we have
\[
d\om{-1}\;\equiv\; 0 \qquad(\mbox{mod.\ } \om{-1}
\]
and for $p=0,1,\ldots,k-1$
\[
d\om{p} + [\om{-1}\wedge\om{p-1}]_0 \;\equiv\; 0 \qquad (\mbox{mod.\ }
\om0,\om1,\ldots,\om{p})\;.
\]
Therefore the theorem follows from Proposition~\ref{pro:5.1}. \qed

\begin{example}\label{exa:6.1} We define the 6-dimensional Lie algebra
$\jj=(X_1,X_2,X_3,X_4,X_5,X_6)$ by the following bracket operations:
\begin{eqnarray*}
&& [X_1,X_4]=-X_2\;,\; [X_1,X_5]=-X_1\;,\; 
 [X_2,X_4]=-X_3\;,\;[X_2,X_5]=-X_2\\
&& [X_2,X_6]=-X_2\;,\; [X_3,X_5]=-X_3\;,\; 
 [X_3,X_6]=-2X_3\;,\;[X_4,X_6]=-X_4
\end{eqnarray*}
and otherwise $[X_i,X_j]=0$.  Set $\mm=(X_1,X_2,X_3,X_4)$ and
$\jj_0=(X_5,X_6)$.  Then $(\jj,\mm,\jj_0)$ is a differential Lie algebra
of order 1.  Let $G$ be a connected Lie group corresponding to the
Lie algebra $\jj$ and let $\omega$ be the Maurer-Cartan form of $G$.
Setting $\omega=\sum_{i=1}^6 \om{i}X_i$ we have
\[
\left\{\begin{array}{l}
d\om1 = \om1\wedge\om5\;,\\
d\om2 = \om1\wedge\om4+\om2\wedge(\om5+\om6) \;,\\
d\om3 = \om2\wedge\om4+\om3\wedge(\om5+2\om6) \;,\\
d\om4 = \om4\wedge\om6 \;,\\
d\om5 = 0 \;,\\
d\om6 = 0 \;.
\end{array}\right.
\]
\end{example}
The differential equation \reff{eq:3.1} considered in Section~\ref{sec:3} belongs
to this homogeneous case.  The differential equation is left invariant by the Lie group,
which can be considered as a subgroup of the contact transformation group.

\begin{theorem}\label{the:6.2}
If $\jj=\jj_{-1}+\jj_0+\jj_1$ (direct sum) is a semi-simple graded Lie algebra,
i.e.~$[\jj_i,\jj_j]\subset\jj_{i+j}$ ($i,j=0,\pm 1,\pm 2,\ldots$), where we put
$\jj_p=(0)$ for $p\le -2$ and $p\ge 2$, then $\jj$ has a structure of fundamental
differential Lie algebra.
\end{theorem}

\proof
Let $B$ be the Killinng-Cartan form of $\jj$.  The linear endomorphism $s$
of $\jj$ defined by
\[
s(X_{-1}+X_0+X_1) \;=\; -X_{-1}+X_0-X_1
\qquad (X_{-1}\in\jj_{-1}\,,\, X_0\in\jj_0\,,\,X_1\in\jj_1)
\]
is an involutive automorphism of $\jj$.  Hence
\[
B(X_1,X_0) \;=\; B\bigl(s(X_1),s(X_0)\bigr) \;=\;
B(-X_1,X_0) \;=\; -B(X_1,X_0) \qquad (X_0\in\jj_0\,,\,X_1\in\jj_1)\;.
\]
Therefore we have
\begin{equation}
\label{eq:6.1}
B(\jj_1,\jj_0) \;=\;0\;.
\end{equation}
Let $X_1\in\jj_1$ be an element satisfying $[X_1,\jj_{-1}]=(0)$.  
For $Y_{-1}\in\jj_{-1}$ and $Z_i\in\jj_i$ ($i=-1,0,1$), we have
${\rm ad}(X_1)\circ{\rm ad}(Y_{-1})(Z_{-1}) =0$,
\[
{\rm ad}(X_1)\circ{\rm ad}(Y_{-1})(Z_0) \;=\; 
\bigl[X_1,[Y_{-1},Z_0]\bigr] \;=\; 0
\]
and
\[
{\rm ad}(X_1)\circ{\rm ad}(Y_{-1})(Z_1) \;=\; 
-\bigl[Y_1,[Z_1,X_1]\bigr] -\bigl[Z_1,[X_1,Y_{-1}]\bigr] \;=\; 0\;.
\]
Hence
\begin{equation}
\label{eq:6.2}
B(X_1,\jj_{-1}) \;=\;0\;.
\end{equation}
For $Y_{1}\in\jj_{1}$ and $Z_i\in\jj_i$ ($i=-1,0,1$), we have
\[
{\rm ad}(X_1)\circ{\rm ad}(Y_{1})(Z_{-1}) \;=\; 
-\bigl[Y_1,[Z_{-1},X_1]\bigr] -\bigl[Z_{-1},[X_1,Y_{1}]\bigr] \;=\; 0\;,
\]
${\rm ad}(X_1)\circ{\rm ad}(Y_{1})(Z_{0})=0$ and ${\rm ad}(X_1)\circ{\rm ad}(Y_{1})(Z_{1})=0$.
Hence 
\begin{equation}
\label{eq:6.3}
B(X_1,\jj_{0}) \;=\;0\;.
\end{equation}
These three relations \reff{eq:6.1}, \reff{eq:6.2} and \reff{eq:6.3} yield
$B(X_1,\jj) =0$.  
Since $\jj$ is semi-simple, this implies $X_1=0$.  Therefore
$\jj_1$ can be considered a subspace of ${\rm Hom}(\jj_{-1},\jj_0)$
by the mapping $h:\jj_1\to {\rm Hom}(\jj_{-1},\jj_0)$ defined by
$h(X_1)(X_{-1}) = [X_1,X_{-1}]$, $X_1\in\jj_1$, $X_{-1}\in\jj_{-1}$.
Next, we define the element 
$s\in{\rm Hom}\bigl(\jj,\jj_1\bigl(W_1^0(\jj_1)\bigr)\bigr)$ by the formulas:

\begin{eqnarray*}
&& S(X_{-1})=0 \quad,\quad S(X_1)=0 \quad,\quad
S(X_0)(Y_{-1})=[X_0,Y_{-1}]\quad,\\
&&S(X_0)(Y_0)=\frac12[X_0,Y_0] \quad,\quad 
S(X_0)(Y_{1})=[X_0,Y_{1}]\quad,
\end{eqnarray*}
$X_{-1}, Y_{-1}\in\jj_{-1}$, $X_{0}, Y_{0}\in\jj_{0}$,
$X_{1}, Y_{1}\in\jj_{1}$.  Then we have
\[
[X,Y]\;=\;[X,Y]_0 + (\partial S)(X\wedge Y)
\]
for all $X,Y\in\jj$.\qed

The simple real Lie algebras having the structure stated in Theorem~\ref{the:6.2}
are classified in S.~Kobayashi and T.~Nagano~\cite{koba}.
Among these simple Lie algebras ${\rm sl}(2:\real)$ is the simplest example.
Moreover ${\rm sl}(2:\real)$ has the structure of a differential Lie algebra
which is not fundamental. Set
\[
\mm\;=\;\Biggl\{
\left(\begin{array}{cc}
0 & a \\ b & 0
\end{array}\right)
 \Biggm| a,b\in\real\Biggr\}
 \quad,\quad
 \jj_o\;=\;\Biggl\{
\left(\begin{array}{cc}
u & 0 \\ 0 & -u
\end{array}\right)
 \Biggm| u\in\real\Biggr\}\quad.
 \]
$({\rm sl}(2:\real),\mm,\jj_0)$ is a differential Lie algebra of order 1.  Let $\omega$ be
the Maurer-Cartan form of ${\rm SL}(2:\real)$ and set
\[
\omega\;=\; \om1 \left(\begin{array}{cc}
0 & 1\\ 0 & 0
\end{array}\right) +
\om2 \left(\begin{array}{cc}
0 & 0\\ 1 & 0
\end{array}\right) -
\om3 \left(\begin{array}{cc}
0 & 1\\ 0 & -1
\end{array}\right)\;.
\]
We have the structure equation
\[
\left\{\begin{array}{l}
d\om1 = -2\om1\wedge\om3\;,\\
d\om2 = 2\om2\wedge\om3 \;,\\
d\om3 = \om1\wedge\om2  \;.
\end{array}\right.
\]
For example, we can choose $\om1,\om2,\om3$ as follows:
\[
\left\{\begin{array}{l}
\om1 = {\rm e}^{2z} \,dx\;,\\
\om2 = {\rm e}^{-2z} \bigl(dy-\bigl(y^2+a(x)y+b(x)\bigr)dx\bigr)\;,\\
\om3 =dz-\bigl(y+\frac12 a(x)\bigr)dx \;,
\end{array}\right.
\]
where $a(x)$ and $b(x)$ denote two arbitrary differentiable functions of the
variable $x$.  Therefore we can see that ${\rm SL}(2:\real)$ corresponds
to the family of ordinary differential equations of Riccati type.

\section{Cartan example}\label{sec:7}
Now we can consider the involutive system of second-order differential equations
which is one of the typical examples in Cartn's paper~\cite{cartpfaff}:
\begin{equation}
  \label{eq:7.1}
  {\partial^2 z\over\partial x^2}=\frac13
  \Bigl({\partial^2 z\over\partial y^2}\Bigr)^3 \quad , \quad
{\partial^2 z\over \partial x\,\partial y} = \frac12 \Bigl({\partial^2 z\over\partial y^2}\Bigr)^2\quad.
\end{equation}

Setting on $\real^6=\{(x,y,z,p,q,t)\}$
$\om1=dx$, $\om2=dy+t\,dx$, $\om3=dz-p\,dx-q\,dy$,
$\om4=dp-t\,dq+\frac16 t^3\,dx+\frac12 t^2\,dy$, $\om5=dq-\frac12 t^2\,dx-t\,dy$
and $\om6=dt$, we have the structure equations of the system:
\begin{equation}
\label{eq:7.2}
\left\{\begin{array}{l}
d\om1 = 0\;,\\
d\om2=-\om1\wedge\om6\;,\\
d\om3 = \om1\wedge\om4+\om2\wedge\om5 \;,\\
d\om4 = \om5\wedge\om6 \;,\\
d\om5 = \om2\wedge\om6 \;,\\
d\om6 = 0\;.
\end{array}\right.
\end{equation}
which is of type $W_2^0$ inn Example~\ref{exa:5.2} (2)(ii)(c).  Th integration 
of the system is deeply related to the structure equation~\reff{eq:7.2}.
It is easy to see that $(\om6,\om2,\om5,\om4,\om3)$ forms a solvable
system of the Pfaffian system $\Sigma=\langle\om3,\om4,\om5\rangle$.  
Therefore we can obtain by quadrature five independent first integrals
of the Cauchy characteristic system of $\Sigma$:
\begin{eqnarray*}
&&u_1=z-xp+xqt+\frac16 x^2t^3 \quad,\quad 
u_2=p-qt+\frac16 xt^3+\frac12 y t^2\\
&& u_3=q-\frac12 y t^2 \quad,\quad u_4=y+xt
\quad,\quad u_5=t\quad,
\end{eqnarray*}
so that the system $\Sigma$ is expressed as
\[
\left\{\begin{array}{rcl}
\om3-x\om4 &=& du_1-(u_3+u_4\,u_5)\,du_4\;,\\
\om4 &=& du_2 +u_3\,du_5\;,\\
\om5 &=& du_3+u_4\,du_5\;.
\end{array}\right.
\]
By this expression, we can construct the general integral surfaces of the given system
(see~\cite[\S 38]{cartpfaff}).

Let $G$ be the Lie subgroup of ${\rm GL}(6:\real)$ consisting of matrices of the form
\[
\left(
\begin{array}{cccccc}
* & * & * & * & * & 0\\
{*} & * & * & * & * & 0\\
0 & 0 & * & 0 & 0 & 0\\
0 & 0 & * & * & * & 0\\
0 & 0 & * & * & * & 0\\
{*} & * & * & * & * & *
\end{array}
\right)\;,
\]
where $*$ is an element of $\real$.  Let $B_G$ be the $G$-structure defined by
the dual frame of $(\om1,\om2,\ldots,\om6)$.  Then a contact transformation leaving
the equations~\reff{eq:7.2} invariant induces an automorphism of this $G$-structure
and vice versa.  The structure group $G$ can be reduced to the Lie subgroup
$G_7$ whose Lie algebra $\jj_7$ is given as follows:
\[
\jj_7\;=\;
\left\{
\left(
\begin{array}{cccccc}
a_1-a_2 & -a_3 & a_4 &0 &0 & 0\\
0&a_1 &a_5 & 0& \frac43 a_3& 0\\
0& 0& 2a_1+a_2& 0& 0& 0\\
0& 0& a_6 &a_1+2a_2 & 0& 0\\
0& 0& a_7& a_3& a_1+a_2& 0\\
0& -a_6& 0& -a_5& \frac43 a_7& a_2
\end{array}
\right)\;,\; a_i\in\real\,,\, i=1,2,\ldots,7
\right\}\;.
\]
The usual prolongations of $\jj_7$ satisfy $\dim\jj_7^{(1)}=1$ and $\jj_7^{(P)}=\{0\}$
($p\ge 2$).  Therefore the local automorphism group of the given system
with respect to the group $G$ is of finite type.

The structure equation~\reff{eq:7.2} determines the Lie algebra
$\mm=(X_1,X_2,X_3,X_4,X_5,X_6)$ with the bracket operation
\begin{eqnarray*}
&& [X_1,X_4]=-X_3\;,\; [X_1,X_6]=-X_2\;,\; 
 [X_2,X_5]=-X_3\\
&& [X_2,X_6]=-X_5\;,\; [X_5,X_6]=-X_4
\end{eqnarray*}
and otherwise $[X_i,X_j]=0$.  This Lie algebra has a structure of fundamental
graded Lie algebra in the sense of N.~Tanaka~\cite{tana}.  We set
$\jj_{-1}=(X_1,X_6)$, $\jj_{-2}=(X_2)$, $\jj_{-3}=(X_5)$, $\jj_{-4}=(X_4)$,
$\jj_{-5}=(X_3)$.  Then $\mm=\jj_{-5}+\jj_{-4}+\jj_{-3}+\jj_{-2}+\jj_{-1}$ 
(direct sum) is a fundamental graded algebra of the 5th kind.  The structure of
the automorphism group with respect to the group $G$ is given by
Tanaka's prolongation method~\cite{tana}.  We can verify that the graded Lie algebra
$\jj$ prolonged from $\mm$ with respect to the Lie algebra of $G$ has the following
structure:
\begin{itemize}
\item[(1)] $\jj=\sum_{p=-5}^5 \jj_p$ (direct sum);
\item[(2)] $\jj_0$ is a Cartan subalgebra of $\jj$;
\item[(3)] $\dim\jj_0=2$, $\dim\jj_{\pm1}=2$ and $\dim\jj_{\pm p}=p$
for $p=2,3,4,5$;
\item[(4)] $\jj$ is isomorphic to the exceptional simple real Lie algebra of 
dimension 14.
\end{itemize}
In order to transform the involutive systems which admit this 14-dimensional
simple Lie group of contact transformations to the standard form~\reff{eq:7.1},
we need the integration of a system of differential equations associated with 
the simple group.

\bibliographystyle{plain}
\nocite{*}
\bibliography{mabiblio}
\end{document}